\documentclass{article}
\usepackage[utf8]{inputenc}
\usepackage[]{overpic}
\usepackage{lineno}
\usepackage{bm}
\usepackage{amssymb}
\usepackage{enumerate}
\usepackage{amsmath}
\usepackage{amsthm}
\usepackage{epsfig}
\usepackage{graphicx}
\usepackage{graphics}
\usepackage{float}
\usepackage{subfigure}
\usepackage{multirow}
\usepackage{color}
\usepackage{lineno}
\usepackage{fullpage}
\usepackage[normalem]{ulem} 
\usepackage{makeidx}
\usepackage{xspace}
\usepackage{wrapfig}
\usepackage{mathtools}
\usepackage{hyperref}
\usepackage[table,xcdraw]{xcolor}

\usepackage{booktabs}
\usepackage{csvsimple}

\usepackage{algorithm}
\usepackage{algorithmic}
\usepackage{pgfplotstable,filecontents}
\pgfplotsset{compat=1.9}

\newcommand {\mm}[1] {\ifmmode{#1}\else{\mbox{\(#1\)}}\fi}

\newcommand{\R}{\mathbb{R}}
\newcommand{\norm}[1]{\left\lVert#1\right\rVert}

\DeclarePairedDelimiterX\set[1]\lbrace\rbrace{#1}
\makeindex

\title{Sparsifying Priors for Bayesian Uncertainty Quantification\\ in Model Discovery}%
\author{Seth~M. Hirsh$^*$, %
  David~A. Barajas-Solano$^\dag$, %
  J.~Nathan Kutz$^\ddag$\\[.1in]
  {{$^*$ Department of Physics, University of Washington, Seattle, WA}}\\
  {$^{\dag}${ Pacific Northwest National Laboratory, Richland, WA}}\\
{$^{\ddag}$ { Department of Applied Mathematics, University of Washington, Seattle, WA}} }%
\date{\today}

\begin{document}

\maketitle

\begin{abstract}
We propose a probabilistic model discovery method for identifying ordinary differential equations (ODEs) governing the dynamics of observed multivariate data.
Our method is based on the {\em sparse identification of nonlinear dynamics} (SINDy) framework, in which target ODE models are expressed as a sparse linear combinations of pre-specified candidate functions.
Promoting parsimony through sparsity in SINDy leads to interpretable models that generalize to unknown data.
Instead of targeting point estimates of the SINDy (linear combination) coefficients, in this work we estimate these coefficients via sparse Bayesian inference.
The resulting method, {\em uncertainty quantification SINDy} (UQ-SINDy), quantifies not only the uncertainty in the values of the SINDy coefficients due to observation errors and limited data, but also the probability of inclusion of each candidate function in the linear combination.
UQ-SINDy promotes robustness against observation noise and limited data, interpretability (in terms of model selection and inclusion probabilities), and generalization capacity for out-of-sample forecast.
Sparse inference for UQ-SINDy employs Markov Chain Monte Carlo, and we explore two sparsifying priors: the {\em spike-and-slab prior}, and the {\em regularized horseshoe prior}.
We apply UQ-SINDy to synthetic nonlinear data sets from a Lotka-Volterra model and a nonlinear oscillator, and to a real-world data set of lynx and hare populations.
We find that UQ-SINDy is able to discover accurate and meaningful models even in the presence of noise and limited data samples.

\end{abstract}

\section{Introduction}
In recent years there has been a rapid increase in measurements gathered from complex nonlinear dynamics for which their governing equations are unknown.
A key challenge is to discover explicit representations of these equations, which can then be used for system identification, forecasting and control.
Measurements are often compromised by noise or may exhibit chaotic behavior, in which case it is critical to quantify how uncertainty affects the model discovery process.
To address this challenge, we introduce the {\em uncertainty quantification sparse identification of nonlinear dynamics} (UQ-SINDy) framework, which leverages sparsity promotion in a Bayesian probabilistic setting to extract a parsimonious set of governing equations.
Our method provides uncertainty estimates of both the parameter values and the inclusion probabilities for different terms in the models. 
    
Discovery of governing equations plays a fundamental role in the development of physical theories.
With increasing computing power and data availability in recent years, there have been substantial efforts to identify the governing equations directly from data~\cite{bongard2007automated,schmidt2009distilling,yang2020bayesian}.
There has been particular emphasis on parsimonious representations because they have the benefits of promoting interpretibility and generalizing well to unknown data~\cite{bai2015low,brunton2014compressive,brunton2016discovering,mackey2014compressive,ozolicnvs2013compressed,proctor2014exploiting,tran2017exact,wang2011predicting}.  The SINDy method was propoosed in~\cite{brunton2016discovering}, which leverages dictionary learning and sparse regression to model dynamical systems. This approach has been successful in modeling a diversity of applications, including in chemistry~\cite{hoffmann2019reactive}, optics~\cite{sorokina2016sparse}, engineered systems~\cite{li2019discovering},  epidemiology~\cite{horrocks2020algorithmic}, and plasma physics~\cite{dam2017sparse}. Furthermore, there has been a variety of modifications, including improved robustness to noise~\cite{champion2019data, champion2020unified,kaheman2020automatic}, generalizations to partial differential equations~\cite{raissi2018hidden,rudy2017data,rudy2019data}, multi-scale physics~\cite{champion2019discovery}, and libraries of rational functions~\cite{mangan2016inferring,kaheman2020sindy}.

Although these methods identify the equations, measurements often contain observation errors, which may imperil the predictive capacity of learned models.
A common approach to remedy this is to use the Bayesian probability framework where uncertainty is quantified in terms of probability and where priors are employed to encode assumptions and knowledge about model parameters~\cite{gelman2013bayesian,west2006bayesian}.
Bayesian methods have been widely used for uncertainty quantification in time series models, with applications to weather forecasting~\cite{abramson1996hailfinder,elsner2004hierarchical,yu2013bayesian}, disease modeling~\cite{best2005comparison,lawson2013bayesian,yuen2002bayesian}, traffic flow~\cite{castillo2008predicting,sun2006bayesian,zheng2006short}, and finance~\cite{gerlach2011bayesian,ticknor2013bayesian,wright2008bayesian}, among many others.
More recently, these methods have been incorporated into model discovery frameworks, exhibiting state-of-the-art performance for system identification in the presence of noise~\cite{galioto2020bayesian,niven2020bayesian,yang2020bayesian}.
Although these methods provide a range of possible values, realizations of these models are in general not sparse and consequently lack the capability to identify relevant terms in the model.

Sparse regression is a popular tool to identify a small subset of variables that explain the data.
However, finding the true minimum is computationally intractable in practice.
In the frequentist setting, a popular solution is to use the Lasso, which corresponds to an $l_1$ penalty term~\cite{tibshirani1996regression}.
In the Bayesian setting, sparsity is generated by fundamentally different mechanisms.
Most notably, although the corresponding prior (the Laplace prior) shares the same maximum likelihood estimator as the Lasso~\cite{park2008bayesian}, the distribution has fat tails and thus does not produce sparse realizations~\cite{castillo2015bayesian}.
The spike and slab model remedies this by explicitly using Bernoulli variables to determine whether a term is present in the model, and has become the leading method for incorporating sparsity in the Bayesian framework~\cite{mitchell1988bayesian,ishwaran2005spike,madigan1994model}.
One disadvantage to this prior however is its dependence on discrete variables, which makes inference prohibitively expensive for high-dimensional systems.
Smooth approximations, such as the horseshoe~\cite{carvalho2009handling,carvalho2010horseshoe}, Horseshoe+~\cite{bhadra2017horseshoe+}, regularized horseshoe~\cite{piironen2017sparsity}, Dirichlet-Laplace~\cite{bhattacharya2015dirichlet}, and R2-D2 priors~\cite{zhang2016high}, have been shown to yield performance comparable to the spike and slab model.
For this work we will primarily focus on the regularized horseshoe prior, also known as the Finnish horseshoe.

In this work, we propose the UQ-SINDy framework, which provides uncertainty estimates of both the parameter value and inclusion probabilities and promotes sparsity in realizations of the model.
This model leverages advances sparsity and Bayesian approaches for solving ODEs to achieve this goal.
In Sections~\ref{subsec:SINDy} and~\ref{subsec:bayesian_inference} we review the SINDy method and Bayesian inference for ordinary differential equations, respectively.
In Section~\ref{subsec:sparse_priors} we review sparsity promoting priors, namely the spike and slab and regularized horseshoe priors, and compare their performance to the Laplace prior.
In Section~\ref{subsec:method}, we introduce two sparsity promoting Bayesian methods, spike and slab SINDy and regularized horseshoe SINDy.
In Sections~\ref{subsec:synthetic_example} and~\ref{subsec:lynx_hare}, we illustrate these methods on two synthetic nonlinear data sets, a Lotka Volterra model and nonlinear oscillator, and one real-world example of lynx and hare population data.
We find that these methods are able to extract accurate and meaningful Bayesian models even in the presence of significant noise and sparse samples.
These results are summarized and future improvements are discussed in Section~\ref{sec:conclusions}.

\section{Background}

\begin{figure}
\includegraphics[width=1\textwidth]{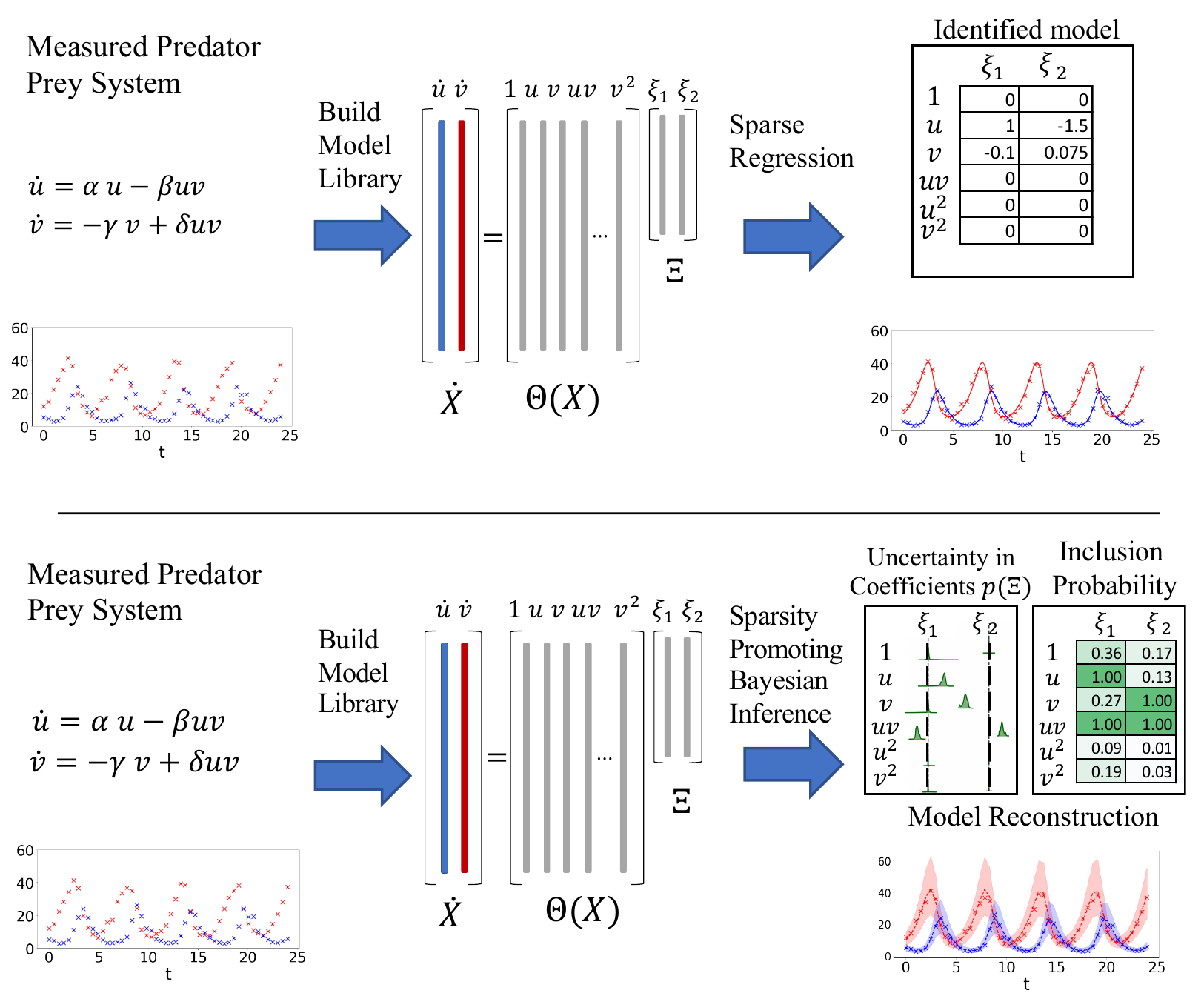}
\vspace*{-.2in}
\caption{Comparison of SINDy algorithm and UQ-SINDy. (Top) Schematic of SINDy algorithm. A dynamical system governed by unknown governing equations is measured. Next, we computed the derivative of the time series $\dot{\bm{X}}$ and construct a library $\bm{\Theta(X)}$ of candidate terms. Last, we perform sparse regression to identify the terms in the library that best explain the time series.
  (Bottom) Schematic of UQ-SINDy algorithm. A dynamical system governed by unknown governing equations is measured. Next, we posit a SINDy library $\bm{\Theta(X)}$ of candidate terms. Last, we perform sparsity promoting Bayesian inference to compute the inclusion probability and the posterior distribution of each term in the SINDy library. An ensemble of reconstructions can then be computed, which quantify the credibility of predictions.
}
\label{fig:sindy}
\end{figure}

The UQ-SINDy framework is based on several recent developments in the fields of sparse regression, ordinary differential equations, and Bayesian inference, and we review these contributions here. In Section \ref{subsec:SINDy} we introduce the SINDy algorithm, which employs sparse regression to identify governing equations in the frequentist setting. In Section~\ref{subsec:bayesian_inference}, we review Bayesian inference for ordinary differential equations. In Section~\ref{subsec:sparse_priors}, we review three different priors for sparse inference---the Laplace, spike and slab, and regularized horseshoe priors---and compare their benefits and drawbacks.

\subsection{Sparse identification of nonlinear dynamics (SINDy)}\label{subsec:SINDy}

The SINDy method is a recently developed technique that leverages sparse regression to identify the governing equations from a given time series (Figure~\ref{fig:sindy}).
We consider a system with state $\bm{\bm{x}}(t) = [x_1(t), x_2(t) , \ldots x_d (t)]^{\top} \in \R^d$ governed by the differential equation
\begin{equation*}
    \dot{\bm{x}} = \bm{f}(\bm{x}),
\end{equation*}
for some unknown function $\bm{f} \colon \R^d \to \R^d$.
The system's state is observed at the discrete times $t = t_1, \ldots, t_n$.
The goal of SINDy is to discover $\bm{f}$ from these observations.

To do so, we postulate that $\bm{f}$ can be written as a linear combination of a library of $l$ candidate functions $\theta_i \colon \mathbb{R}^d \to \mathbb{R}$, $i \in [1, l]$.
For example, a commonly used library is the polynomial library
\begin{equation*}
  \bm{\Theta}(\bm{x}) =
  \begin{bmatrix}
    x_1(t) & x_2(t) & x_1^2(t) & x_1(t) x_2(t) & \cdots
  \end{bmatrix} \in \mathbb{R}^l.
\end{equation*}
Next, we define $\bm{X} = [\bm{x}(t_1), \bm{x}(t_2), \ldots,  \bm{x}(t_n)]^{\top} \in \R^{n \times d}$ as the collection of observed state snapshots, and also define the matrix of library terms evaluated at the observation times,
\begin{equation*}
  \bm{\Theta}(\bm{X}) = \begin{bmatrix}
    \bm{\Theta}(\bm{x}(t_1))^{\top} & \bm{\Theta}(\bm{x}(t_2))^{\top} & \cdots & \bm{\Theta}(\bm{x}(t_n))^{\top}
  \end{bmatrix}^{\top} \in \R^{n \times l}.
\end{equation*}
We then measure or compute the time derivative of the data $\bm{X}$ and solve the following equation for $\bm{\Xi} \in \R^{l \times d}$,
\begin{equation}
  \dot{\bm{X}} = \bm{\Theta}(\bm{X}) \bm{\Xi},
  \label{eq:SINDy_problem}
\end{equation}
where $\bm{\Xi}$ denotes the matrix of linear combination coefficients, or {\em SINDy coefficients}.
A key assumption of SINDy is that $\bm{f}$ may be represented by a small number of library terms, so that the matrix $\bm{\Xi}$ is sparse.
Thus, \eqref{eq:SINDy_problem} is typically solved through sparse regression, using minimization techniques such as sequential least squares thresholding (STLSQ)~\cite{brunton2016discovering}, Lasso~\cite{tibshirani1996regression}, or a relaxed formulation~\cite{champion2020unified}. The SINDy procedure yields the set of identified nonlinear differential equations
\begin{equation}
    \dot{\bm{x}}^{\top} = \bm{\Theta} (\bm{x}) \bm{\Xi}.
\end{equation}
Once identified, this system of differential equations may be used for system identification, prediction, and control.

\subsection{Bayesian inference for data-driven discovery}\label{subsec:bayesian_inference}

Suppose we have the data set $(\bm{X}, \bm{y})$ for which we would like to fit the linear regression model
\begin{equation}
  \bm{y} = \bm{\beta}^{\top} \bm{X} + \bm{\epsilon}
  \label{eq:regression_problem}
\end{equation}
where $\bm{\epsilon} \sim \mathcal{N}(0, \sigma^2 \bm{I})$ is a vector of independent, identically distributed Gaussian measurement noise with unknown standard deviation $\sigma$.
In the Bayesian setting, our goal is to determine the posterior distribution of $\bm{\beta}$ and $\sigma$ conditioned on the data, i.e. $p(\bm{\beta}, \sigma | \bm{X}, \bm{y})$.
To compute this distribution, we leverage Bayes' rule,
\begin{equation*}
  p(\bm{\beta}, \sigma | \bm{X}, \bm{y}) \propto {p(\bm{y} | \bm{\beta}, \bm{X}) \, p(\sigma) \, p(\bm{\beta})},
\end{equation*}
where $p(\bm{y} | \bm{\beta}, \bm{X})$ denotes the data likelihood, and $p(\sigma$) and $p(\bm{\beta})$ denote the prior distribution of the noise standard deviation and the regression coefficients.
These prior distributions incorporate any available domain knowledge about the distribution of the noise standard deviation and the $\beta_j$s.

In this work we are interested in identifying ODE models from noisy data. In particular, given noisy time series and a SINDy model of the form $\dot{\bm{x}}^{\top} = \bm{\Theta}(\bm{x}) \bm{\Xi}$, our goal is to compute the posterior distribution of the initial conditions $\bm{x}_0$ and SINDy coefficients $\bm{\Xi}$.
We assume that the data set $\bm{X}$ consists of $n$ noisy snapshots of the observed dynamics, that is $\bm{X} = [\bm{y}_1, \bm{y}_2, \ldots, \bm{y}_n]^{\top} \in \mathbb{R}^{n \times d}$, where $\bm{y}_i \in \mathbb{R}^d$ is the noisy snapshot of the system state at time $t = t_i$.
For a given probabilistic model of the observation noise, the data is modeled as deviations from the SINDy predictions; for example, for additive noise models,
\begin{equation}
  \bm{y}_i^{\top} = \bm{x}^{\top}_0 + \int_0^{t_i} \bm{\Theta}(\bm{x}(t')) \bm{\Xi} \, dt' + \bm{\epsilon}_i,
  \label{eq:ode_problem}
\end{equation}
where $\bm{\epsilon}_i$ denotes the additive noise for the $i$th snapshot.
Bayes' rule then takes the form
\begin{equation}
  p(\bm{\Xi}, \bm{x}_0, \bm{\phi} | \bm{X}) \propto p(\bm{X}| \bm{\Xi}, \bm{x}_0, \bm{\phi}) \, p(\bm{\phi}) \, p(\bm{\Xi}) \, p(\bm{x}_0)
  \label{eq:bayes},
\end{equation}
where $\bm{\phi}$ denotes auxiliary variables of the probabilistic model such as the noise standard deviation.
The data likelihood $p(\bm{X}| \bm{\Xi}, \bm{x}_0, \bm{\phi})$ is given by the chosen observation model (e.g., by \eqref{eq:ode_problem} and the distribution of the noise for additive observation noise).

Computing the posterior distribution~\eqref{eq:bayes} is in general not analytically tractable, in which case sampling-based methods such as Markov Chain Monte Carlo (MCMC) may be used.
Once the posterior distribution has been approximated, we may then compute state reconstructions and forecasts conditioned on the observed data~\cite{gelman1996posterior,tran2016edward}. Specifically, to estimate the distribution of predicted values of $\bm{x}$ at an arbitrary time $t$, we marginalize the data likelihood times the posterior distribution over $\bm{\Xi}$, $\bm{x}_0$, and $\bm{\phi}$, that is,
\begin{equation}
    p(\bm{x}(t) | \bm{X}) = \int  p(\bm{x} (t) | \bm{\Xi}, \bm{x}_0, \bm{\phi}) \, p(\bm{\Xi}, \bm{x}_0, \bm{\phi} | \bm{X}) \, d\bm{\Xi} \, d\bm{x}_0 \, d\bm{\phi},
\label{eq:ppc}
\end{equation}
The distribution $p(\bm{x}(t) | \bm{X})$ is referred to as the {\em posterior predictive distribution (PPD)}.
The integral in~\eqref{eq:ppc} can be approximated via sampling by taking the expectation of the data likelihood over posterior samples drawn via MCMC.

\subsection{Sparsity promoting priors}\label{subsec:sparse_priors}

\begin{figure}
\includegraphics[width=\textwidth]{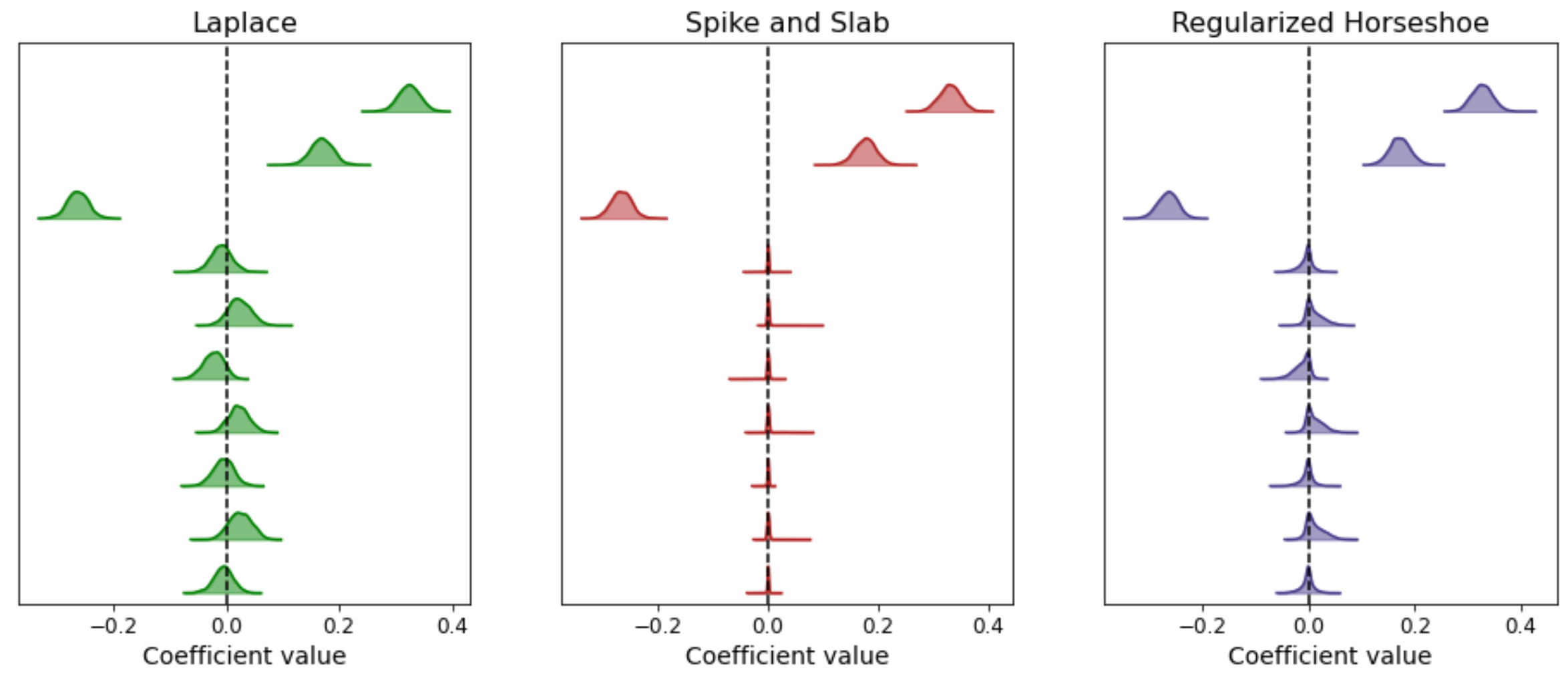}
\caption{Comparison of posterior distributions for Laplace, spike and slab, and regularized horseshoe priors for a linear regression problem. Both the spike and slab and regularized horeshoe priors promote sparsity in the posterior distributions, while the Laplace prior does not. 
}
\label{fig:priors}
\end{figure}

Consider the regression problem in~\eqref{eq:regression_problem}. In many cases, we assume only a few components of $\bm{x}_i$ are relevant for predicting $y_i$, in which case we expect $\bm{\beta}$ to be sparse. In the Bayesian setting, multiple sparsity-inducing priors have been proposed.
We describe a few of these approaches below, namely the Laplace, spike and slab and regularized horseshoe priors.

\subsubsection{Laplace prior}
Originally proposed by Laplace~\cite{de1774memoire}, the Laplace distribution, also known as the double exponential distribution~\cite{gelman2013bayesian}, corresponds to the probability distribution function (PDF) $f(x | \mu, b)$ given by
\begin{equation*}
  f(x | \mu, b) = \frac{1}{2b} \exp \left \{ - \left | \frac{x - \mu}{b} \right | \right \}.
\end{equation*}
Most notably, maximum a posterior (MAP) estimation for this prior corresponds to regression with $\ell_1$ regularization, that is,~\cite{park2008bayesian}, 
\begin{equation*}
    \hat{\bm{\beta}}^{\text{Laplace}} = \operatorname{arg\,max}_{\bm{\beta}} \, p(\bm{y} | \bm{\beta}, \bm{X}) \, p(\bm{\beta}) = \operatorname{arg\,min}_{\bm{\beta}} \, \lvert| \bm{y} - \bm{\beta}^{\top} \bm{X} \rvert_2^2 + \lambda \norm{\bm{\beta}}_1.
\end{equation*}
In the frequentist setting, solving this regression problem, known as the LASSO problem, has been shown to yield sparse solutions for $\bm{\beta}$~\cite{tibshirani1996regression}. This sparsifying behavior of the Laplace distribution is attributed to the fact that for values of $x$ smaller than $b$, the distribution is sharply peaked, thus pushing many terms toward $0$. Additionally, for values of $x$ greater than $b$, the distribution has longer tails than the Gaussian distribution, allowing elements to escape significant shrinkage.  

Although $l_1$ regularization induces sparsity in the frequentist case, in the Bayesian setting realizations of the corresponding posterior distributions are not sparse~\cite{castillo2015bayesian}.
In particular, in the Bayesian setting we must consider the whole distribution simultaneously.
With the Laplace prior, every $\beta_j$ has probability mass simultaneously pushed toward and away from the origin, forcing relevant $\beta_j$s to shrink toward the origin and irrelevant terms to have significant probability mass far away from the origin.

To illustrate this we generate 400 data samples ($\bm{x}_i,y_i$) that satisfy~\eqref{eq:regression_problem}, where $\bm{x}_i \sim \mathcal{N}(0,1) \in \R^{10}$, $\epsilon_i \sim \mathcal{N}(0,0.5^2)$, and $\bm{\beta} \in \R^{10}$ chosen to be the sparse vector 
\begin{equation*}
    \bm{\beta} = [0.3, 0.2, -0.3, 0, 0, 0, 0, 0, 0 ,0]^{\top}.
\end{equation*}
We perform Bayesian inference to estimate $\bm{\beta}$ using a Laplace prior, and we plot the resulting posterior distribution in Figure~\ref{fig:priors}. We note that when using the Laplace prior, the posterior distributions are centered about the true value $\bm{\beta}$. However, many distributions are peaked at nonzero values, making it difficult to differentiate between relevant and irrelevant variables. Further, due to the wide widths of of all the distributions, samples from this posterior distribution will not be sparse. To better enforce sparsity in a Bayesian setting and induce sparse realizations, the distribution of each $\beta_j$ must either be fully shrunk towards the origin or pushed away from the origin. In Sections~\ref{subsubsec:spike_slab} and~\ref{subsubsec:regularized_horseshoe} we discuss two priors that satisfy these properties.

\subsubsection{Spike and slab prior}\label{subsubsec:spike_slab}

The spike and slab prior is one of the most popular sparsifying priors and is typically referred to as the ``gold standard''~\cite{mitchell1988bayesian,ishwaran2005spike,madigan1994model} sparsity-inducing prior in the Bayesian setting.
For this prior, each $\beta_j$ is generated using the hierarchical model
\begin{align*}
  \beta_j | \lambda_j &\sim \mathcal{N}(0, c^2) \lambda_j, \\
  \lambda_j &= \operatorname{Ber}(\pi),
  \label{eq:spikeandslab}
\end{align*}
where $\operatorname{Ber}(\pi)$ denotes the Bernoulli distribution with probability of success $\pi$.
Here, $\pi$ is the prior probability that $\lambda_j$ is $1$. Otherwise $\lambda_j$ is $0$. From this it can be seen that if $\lambda_j$ is $1$, then the $j$th term belongs to the model and $\beta_j$ follows the ``slab'' distribution, a normal distribution with variance $c^2$. If $\lambda_j$ is $0$, then the $j$th term is not in the model and $\beta_j$ follows the ``spike'' distribution, a Dirac delta distribution centered at zero.

The distribution may be relaxed to 
\begin{gather*}
    \beta_j | \lambda_j \sim \lambda_j \mathcal{N}(0, c^2) + (1 - \lambda_j) \mathcal{N}(0, \epsilon^2),\\
    \lambda_j = \operatorname{Ber}(\pi),
\end{gather*}
where $\epsilon \ll c$. This is similar to before, except when $\lambda_j = 0$, $\beta_j$ follows a narrow normal distribution with variance $\epsilon^2$.

The spike and slab prior for $\bm{\beta}$ is very intuitive and has shown robust performance in practical applications. In Figure~\ref{fig:priors}, we plot the resulting posterior distribution for the example in Section~\ref{subsubsec:spike_slab}. Most notably we see that similar to the Laplace prior, the spike and slab prior extracts out wide distributions for the three nonzero coefficients. For the seven zero coefficients, on the other hand, the distribution is sharply spiked at the origin. Consequently, any samples drawn from this posterior distribution will be truly sparse. Compared to the Laplace distribution, the nonzero terms are much more easily identifiable. Furthermore, the mean of $\lambda_j$ correspond to the estimate of the ``inclusion probability'', that is the likelihood that a particular $\beta_j$ is relevant to the model.

Although the spike and slab prior has many beneficial properties, one downside is that because of its discrete nature, inference with this prior requires exploring the combinatorial space of possible models. To address this challenge, many smooth approximations to the spike and slab prior distribution have been proposed. We discuss one recent approach in Section~\ref{subsubsec:regularized_horseshoe}.

\subsubsection{Regularized horseshoe prior}\label{subsubsec:regularized_horseshoe}

The horseshoe prior and the recently developed regularized horseshoe prior are smooth priors that have shown comparable performance to the spike and slab model. The horseshoe is defined as the hierarchical prior
\begin{align*}
  \beta_i | \lambda_i, \tau &\sim \mathcal{N}(0, \lambda_i^2 \tau^2), \\
  \lambda_i &\sim \mathrm{C}^{+}(0,1), \\
  \tau &\sim \mathrm{C}^{+}(0,\tau_0),
\end{align*}
where $\rm{C}^{+}(\cdot, \cdot)$ denotes the half-Cauchy distribution~\cite{bhadra2019lasso,carvalho2009handling,carvalho2010horseshoe}. The key intuition behind this prior is that $\tau$ promotes global sparsity, shrinking the posterior distributions of all $\beta_i$s. The $\lambda_i$s, known as the local shrinkage parameters, also have half-Cauchy priors, allowing some of the $\beta_i$s to escape significant shrinkage.
Many analyses have focused on choosing an optimal value for $\tau_0$, and in Piironen~\textit{et al.} values are recommended for sparse linear regression~\cite{piironen2017sparsity}.
In this work we employ $\tau_0 = 0.1$ unless specified otherwise.
We note that decreasing the value of $\tau_0$ increases the sparsity of $\bm{\beta}$ estimates.

One downside of the horseshoe is that relevant terms that ``escape'' shrinkage are not regularized, and thus elements of the posterior distribution may become arbitrarily large. In~\cite{piironen2017sparsity} it was proposed to include a small amount of regularization on each $\lambda_i$, resulting in the {\em regularized horseshoe} prior
\begin{align*}
  \beta_i | \tilde{\lambda}_i, \tau, c &\sim \mathcal{N}(0, \tilde{\lambda}_i^2 \tau^2), \\
  \tilde{\lambda}_i &= \frac{c \lambda_i}{\sqrt{c^2 + \tau^2 \lambda_i^2}} \\
  \lambda_i &\sim \rm{C}^{+}(0,1) \\
  c^2 &\sim \operatorname{Inv-Gamma} \left( \frac{\nu}{2}, \frac{\nu}{2} s^2 \right) \\
  \tau &\sim \rm{C}^{+}(0,\tau_0).
    \label{eq:regularizedhorseshoe}
\end{align*}
where $\operatorname{Inv-Gamma}(\cdot, \cdot)$ denotes the inverse Gamma distribution, and $\nu$ and $s$ are parameters that control the shape of the slab.
For small values of $\lambda_i$, $\lambda_i \tau \ll c$, and $\tilde{\lambda}_i \to \lambda_i$, thus approximating the original horseshoe prior. However, for large values of $\lambda_i$, $\lambda_i \tau \gg c$ and $\tilde{\lambda}_i \to c / \tau$, leading to $\beta_i$ being normally distributed with variance $c^2$. This regularizes $\beta_i$, constraining it to be on the order of $c$.
In this work we employ the values $\nu = 4$ and $s = 2$.

We illustrate the performance of this prior in Figure~\ref{fig:priors} for the example in Section~\ref{subsubsec:spike_slab}. Similar to the spike and slab model, the nonzero coefficients have wide distributions. The zero coefficients, on the other hand are more spiked than those of the Laplace prior, thus resulting in sparser posterior realizations.

Unlike for the spike and slab prior, there is no explicit estimate for the inclusion probabilities. A popular alternative for identifying the relevant terms is to compute the shrinkage factor of the coefficients. 
Specifically, we compute the MAP estimate $\hat{\beta}_i^{Flat}$ with a flat prior (i.e., no prior) and compare it to the MAP estimate with the regularized horseshoe prior, $\hat{\beta}_i^{RH}$. The ratio of these two values is called the shrinkage factor
\begin{equation}
    \kappa_i = \hat{\beta}_i^{RH} / \hat{\beta}_i^{Flat}.
    \label{eq:rh_shrinkage}
\end{equation}
The shrinkage factor of the coefficients has been used to define inclusion ``pseudo-probabilities'' for sparsity-promoting models~\cite{bhadra2019lasso,carvalho2010horseshoe,piironen2017sparsity}.
We employ this approach in this work.
In general these ratios may not lie between $0$ and $1$. 

For our work, we have observed that computing $\hat{\beta}_i^{Flat}$ with flat priors is challenging. To remedy this, we use normal priors $\beta_i \sim \mathcal{N}(0,1)$ instead of flat priors.
Further, we note that~\eqref{eq:rh_shrinkage} can be computed directly from MAP estimates, without having to sample the full posterior distributions. Thus the shrinkage factors can be estimated using optimization techniques instead of full Bayesian inference. However, in practice the associated optimization problems may be nonconvex and highly sensitive to the initial guess. Consequently, for this work we use full Bayesian inference to estimate shrinkage factors.

\section{UQ-SINDY}\label{sec:UQ-SINDy}

In this section, we combine advances in model discovery for dynamical systems and sparsity promoting Bayesian inference to propose the UQ-SINDy framework, which aims to quantify the uncertainty of estimated SINDy coefficients due to measurement, and to estimate the inclusion probabilities for each term in the SINDy library. In particular, within this framwork we introduce two methods: spike and slab SINDy (ss-SINDy) and regularized horseshoe SINDy (rh-SINDy). The ss-SINDy method provides state-of-the art performance for estimating uncertainty of coefficients and inclusion probability, while the rh-SINDy is a smooth approximation that shows comparable performance. We outline this framework below.  

\subsection{Method}\label{subsec:method}

We start with a set of time series measurements $\bm{X} \in \R^{n \times d}$ contaminated by measurement noise.
We assume that our data is governed by the SINDy model
\begin{equation}
  \dot{\bm{x}}^{\top} = \bm{\Theta} ( \bm{x} ) \bm{\Xi}, \quad \bm{x}(0) = \bm{x}_0,
  \label{eq:library}
\end{equation}
for some sparse matrix of SINDy coefficients $\bm{\Xi}$ and initial condition $\bm{x}_0$. Our goal is to determine the posterior distribution $p(\bm{\Xi}, \bm{x}_0, \bm{\phi} | \bm{X})$.
 
\begin{enumerate}[Step 1:]
\item {\bf Construct library}: We posit a library $\bm{\Theta} : \R^{d} \to \R^{l}$ of candidate functions.
We emphasize here that $\bm{\Theta}$ is a symbolic vector function of the system's state $\bm{x}$. This is in constrast to the original SINDy algorithm, in which $\bm{\Theta(X})$ is a fixed matrix computed from the time series data.

Depending on the library, solving the ODE in~\eqref{eq:library} for certain values of initial conditions and parameters may be unstable. Practically, this leads to exploding gradients with respect to SINDy coefficients and initial conditions, and integration steps taken by the ODE solver becoming negligibly small. To remedy this, we add a higher-order polynomial term with a small negative coefficient to the ODE model. For example, for a library of terms up to quadratic order, we add a cubic term, leading to the ODE model
\begin{equation}
  \dot{x}_j = \sum_i \theta_i(\bm{x}) \xi_{i, j} - \varepsilon x^{3}_j,
\end{equation}
where $\xi_{i, j}$ is the $i, j$th element of $\bm{\Xi}$. The parameter $\varepsilon$ is chosen to be sufficiently small so that the ODE is not affected for values of the system's state that lie within the range of the data, but sufficiently large so that $\dot{\bm{x}}$ does not grow too large. In general, if the library $\bm{\Theta}$ includes polynomial terms up to order $n$, we add a term $-\varepsilon x^{n+1}_i$ if $n$ is even, or $-\varepsilon x^{n+2}_i$ if $n$ is odd. This guarantees that the values $\dot{\bm{x}}$ remain finite for both positive and negative values of $\bm{x}$.

\item {\bf Construct model priors and model likelihood.}
  Let $\hat{\bm{x}}(t; \bm{\Xi}, \bm{x}_0)$ denote the SINDy prediction at time $t$ for given values of $\bm{\Xi}$ and $\bm{x}_0$, given by
  \begin{equation*}
    \hat{\bm{x}}^{\top}(t; \bm{\Xi}, \bm{x}_0)  = \bm{x}^{\top}_0 + \int_{t_0}^t \bm{\Theta}(\bm{x}(t')) \bm{\Xi} \, dt'.
  \end{equation*}
  For normally distributed measurement noise, the data likelihood takes the form
  \begin{equation}
    p(\bm{X} | \bm{\Xi}, \bm{x}_0, \bm{\phi}) = \prod^n_{i = 1} \prod^d_{j = 1} \frac{1}{\sigma \sqrt{2 \pi}} \exp \left [ \frac{1}{2 \sigma^2} \left | y_{i,j} - \hat{x}_j(t_i; \bm{\Xi}, \bm{x}_0) \right |^2 \right ].
    \label{eq:normal_likelihood}
  \end{equation}
  For some cases, the values of $\bm{X}$ takes nonnegative values, such as for populations, in which case we may choose to use a lognormal likelihood instead:
  \begin{equation}
    p(\bm{X} | \bm{\Xi}, \bm{x}_0, \bm{\phi}) = \prod^n_{i = 1} \prod^d_{j = 1} \frac{1}{y_{i,j} \sigma \sqrt{2 \pi}} \exp \left [ \frac{1}{2 \sigma^2} \left | \log{y_{i,j}} - \log{\hat{x}_j(t_i; \bm{\Xi}, \bm{x}_0)} \right |^2 \right ].
    \label{eq:lognormal_likelihood}
  \end{equation}

  We must choose priors for the noise level parameter $\sigma$ and the initial conditions $\bm{x}_0$.
  These priors are chosen using knowledge about about the type of parameter (i.e., whether the parameter is nonnegative) and the scales of the data.

\item {\bf Choose a sparsity promoting prior for the SINDy coefficients.}
  Following Section~\ref{subsec:sparse_priors}, for \textit{spike and slab SINDy (ss-SINDy)} we use the hierarchical prior
  \begin{align*}
    \xi_{i, j} | \lambda_j &\sim \mathcal{N}(0, 1) \lambda_{i, j} \alpha_{i, j} \\
    \lambda_{i, j} &\sim \operatorname{Ber}(\pi),
  \end{align*}
  For \textit{regularized horseshoe SINDy (rh-SINDy)} we use the hierarchical prior
  \begin{align*}
    \xi_{i,j} | \tilde{\lambda}_{i,j}, \tau, c &\sim \mathcal{N}(0, 1) \tilde{\lambda}_{i,j} \tau \alpha_{i, j} \\
    \tilde{\lambda}_{i,j} &= \frac{c \lambda_{i,j}}{\sqrt{c^2 + \tau^2 \lambda_{i,j}^2}} \\
    \lambda_{i, j} &\sim \rm{C}^{+}(0,1) \\
    c^2 &\sim \operatorname{Inv-Gamma} \left( \frac{\nu}{2}, \frac{\nu}{2} s^2 \right) \\
    \tau &\sim \rm{C}^{+}(0, \tau_0).
  \end{align*}

  For ss-SINDy, we have that $\bm{\phi}$ consists of the noise level parameter $\sigma$ and the local shrinkage paramaters $\lambda_{i, j}$.
  For rh-SINDy, $\bm{\phi}$ consists of $\sigma$, the $\lambda_{i, j}$s, $c$, and $\tau$.
  The coefficients $\alpha_{i, j}$, which we choose as constants for this analysis, allow us to incorporate any knowledge about the scales of different parameters. For this work we choose $\alpha_i = 1$ unless stated otherwise.

\item {\bf Bayesian Inference.}
  Once the priors and the data likelihood are specified, we employ MCMC to draw samples from the posterior distribution $p(\bm{\Xi}, \bm{x}_0, \bm{\phi} | \bm{X})$.
  Furthermore, we estimate the PPD~\eqref{eq:ppc} for the reconstruction and forecasting tasks of interest. We employ MCMC algorithms as implemented in the Python library PyMC3~\cite{salvatier-2016-pymc3}; specifically, for rh-SINDy we use the No-U-Turn Sampler (NUTS)~\cite{hoffman2014no}, and for ss-SINDy we use the compound step sampler implemented in PyMC3.

  In the UQ-SINDy framework, NUTS leverages the gradients of the SINDy model prediction $\hat{\bm{x}}(t; \bm{\Xi}, \bm{x}_0)$ with respect to $\bm{\Xi}$ and $\bm{x}_0$.
  These gradients are computed using Sunode~\cite{seyboldt-sunode-2021}, a Python wrapper for the CVODES library~\cite{serban2003cvodes} for solving forward and adjoint ODE problems.

\end{enumerate}

\subsection{Examples and applications}\label{subsec:synthetic_example}

In this section, we apply the spike and slab and regularized horseshoe priors in the UQ-SINDy framework and illustrate their performance on three examples: two synthetic data sets and one real-world data set of lynx and hare populations. For each example, we quantify the likelihood of each term of the SINDy library belonging to the underlying dynamical equations, providing both an estimate of the inclusion probability and a distribution of likely values for each SINDy coefficient. We compare these results to the original SINDy algorithm and show that UQ-SINDy significantly outperforms SINDy in identifying the underlying dynamics for noisy observations.

\subsubsection{Lotka-Volterra model}\label{subsec:lotka_volterra}

\begin{figure}
\includegraphics[width=\textwidth]{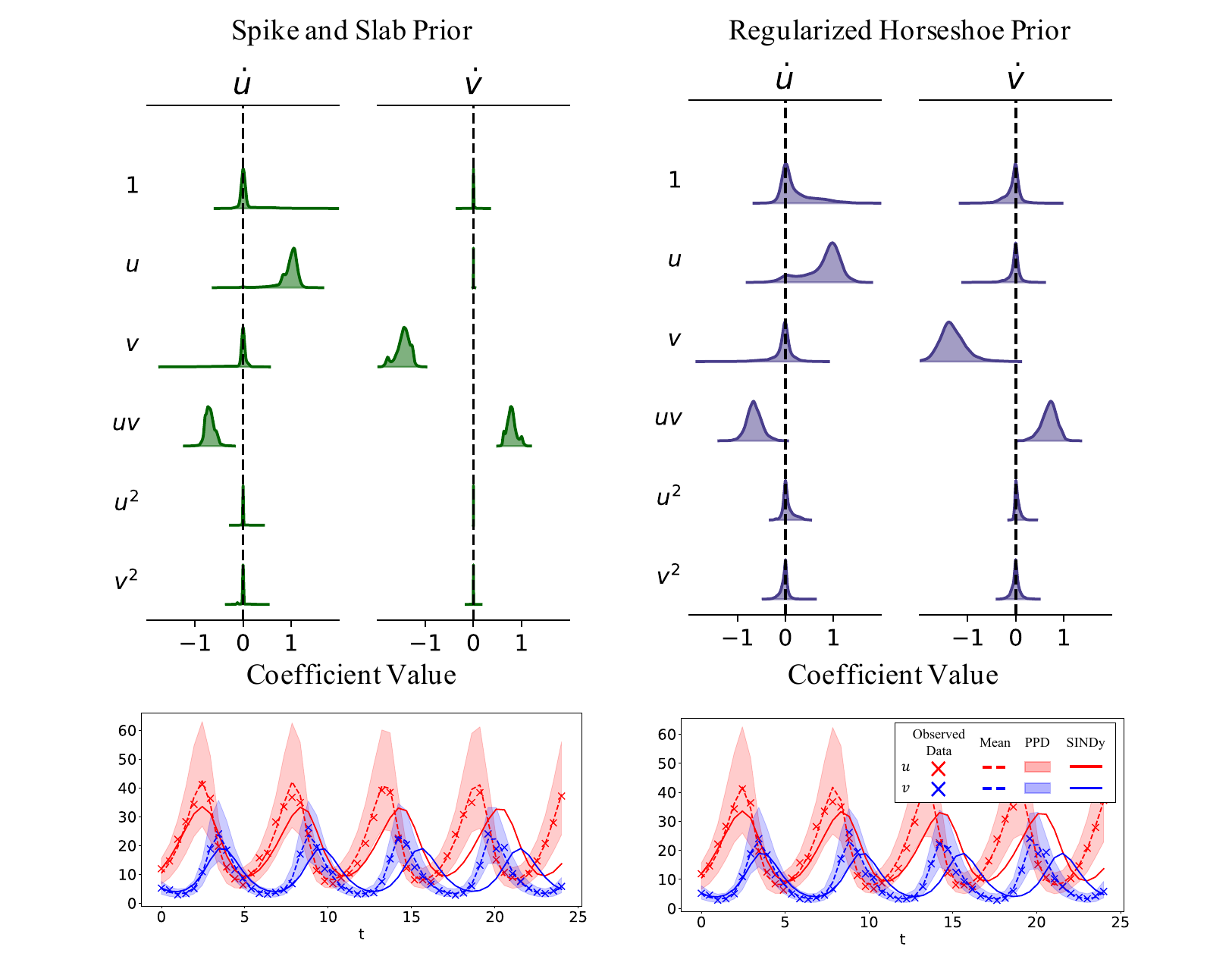}
\caption{UQ-SINDy applied to a synthetic Lotka-Volterra system with lognormal noise.
  (Above) Marginal ss-SINDy and rh-SINDy posterior distributions.
  (Below) Observed (crosses) and predicted time series together with the corresponding PPD means (dashed lines) and 90\% credibility intervals (shaded areas).
  SINDy predictions presented as continuous lines.
}
\label{fig:lotka_volterra}
\end{figure}

We first study data from the Lotka-Volterra model, also commonly refered to as the predator-prey model, which is a popular system used to model the interaction between two competing groups~\cite{goel1971volterra,volterra1927variazioni}. Originally developed by Lotka to model chemical reactions~\cite{lotka2002contribution}, the system has also been studied as a model in economics~\cite{goodwin1982growth} and for biological systems~\cite{kingsland1995modeling,thingstad2000elements,varon1978bacterial}. We explore one real-world example in Section~\ref{subsec:lynx_hare}. 

The Lotka-Volterra model is given by the two nonlinear differential equations
\begin{gather}
\begin{aligned}
    &\dot{u} = \alpha u - \beta u v \\
    &\dot{v} = -\gamma v + \delta u v.
    \label{eq:lotka_volterra}
\end{aligned}
\end{gather}
For this example, we simulate the system with the initial condition $[u_0, v_0] = [10, 5]$ and parameters $\alpha = 1, \beta = 0.1, \gamma = 1.5$, and $\delta = 0.075$, as in~\cite{carpenter2018predator}, which results in a periodic trajectory. We sample 50 snapshots over a time interval of $t \in [0,24]$. Additionally, we contaminate this trajectory with lognormal multiplicative noise with distribution $\text{Lognormal}(0,0.1)$. The lognormal distribution is nonnegative and is commonly used to model observation errors for state variables restricted to with nonnegative values.
The resulting time series is shown in Figure~\ref{fig:lotka_volterra}, from which we see that the trajectory covers approximately four periods of oscillation. 

For this example, we normalize the data as a preprocessing step by dividing each time series (of $x$ and $y$) by the standard deviation of the data. The normalized data is governed by a differential equation of the same form as the unnormalized data, but with modified parameters $\tilde{\alpha} = 1$, $\tilde{\beta} =-0.68$, $\tilde{\gamma} = -1.5$, and $\tilde{\delta} = 0.82$. This preprocessing step can be beneficial for systems in which the parameters are of different orders of magnitude. 

We apply UQ-SINDy for both the spike and slab prior (ss-SINDy) and regularized horseshoe prior (rh-SINDy). We use a library of polynomial terms $\bm{\Theta}(u, v) = [1, u, v, u^2, v^2, u v]$, resulting in a $6 \times 2$ matrix of SINDy coefficients $\bm{\Xi}$.
The SINDy model then reads
\begin{equation*}
  \begin{bmatrix} \dot{u} & \dot{v} \end{bmatrix} = \begin{bmatrix} 1 & u & v & u^2 & v^2 & u v \end{bmatrix} \bm{\Xi}, \quad u(0) = u_0, \quad v(0) = v_0.
\end{equation*}
For the noise level and initial condition we employ the priors $\sigma_u, \sigma_v \sim \operatorname{Lognormal}(\mu = -1, \sigma = 0.1)$ and $u_0, v_0 \sim \operatorname{Lognormal}(\mu = 0, \sigma = 1)$, respectively.

In Table~\ref{table:lotka_volterra}, we present the inclusion probability (for ss-SINDy) and pseudo-probability (for rh-SINDy) of each term in the library. We see significantly higher probabilities for the four true nonzero terms compared to all other terms, indicating that both ss-SINDy and rh-SINDy correctly identify the structure of the governing equation.
We note that although the inclusion pseudo-probabilities are not constrained between zero and one, the relevant terms are easily identified with values near to or greater than 1.

In Figure~\ref{fig:lotka_volterra}, we present the marginal posterior distributions of the SINDy coefficient. From this we immediately see that for both priors, the parameters that belong to the model have broad distributions centered about the true means, while the other 8 terms have narrow peaks centered about 0.
In Table~\ref{table:lotka_volterra}, we compare the posterior modes of the SINDy coefficients against the true values of the model parameters. We additionally apply the original SINDy algorithm to the data. We see that SINDy is unable to identify the correct dynamics due to the presence of observation noise. Furthermore, we note that due to their sparsifying behaviors, the posterior mode of the SINDy coefficients for both the spike and slab and regularized horseshoe priors are close to the true values. 

In Figure~\ref{fig:lotka_volterra}, we present the mean and 90\% credibility interval of the PPDs of the UQ-SINDy reconstructions of the system's states.
Furthermore, we also present the prediction using SINDy (solid lines) and the observed values (crosses). The means of the PPDs for each of the model states are close in value to the true data and provide an accurate continuous reconstruction of the data. In addition, both the regularized horseshoe and spike and slab priors result in similar credibility intervals that bound the true samples. The SINDy reconstruction on the other hand degrades for samples at later times. 

\begin{table}[t]
\centering
\begin{tabular}{|lllll|}
\hline
 & TRUE & SINDy & ss-SINDy & rh-SINDy \\ \hline
$\dot{u} : 1$ & 0 & 0.62 & 0.00 & 0.02 \\
$\dot{v} : 1$ & 0 & 0& 0.00 & -0.01 \\
\rowcolor[HTML]{FFCE93} 
{\color[HTML]{000000} $\dot{u} : u$} & {\color[HTML]{000000} 1} & {\color[HTML]{000000} 0.54} & {\color[HTML]{000000} 1.06} & {\color[HTML]{000000} 0.98} \\
$\dot{v} : u$ & 0 & 0 & 0.00 & 0.00 \\
$\dot{u} : v$ & 0 & -0.49 & 0.00 & -0.01 \\
\rowcolor[HTML]{FFCE93} 
{\color[HTML]{333333} $\dot{v} : v$} & {\color[HTML]{333333} -1.5} & {\color[HTML]{333333} -1.32} & {\color[HTML]{333333} -1.44} & {\color[HTML]{333333} -1.39} \\
\rowcolor[HTML]{FFCE93} 
{\color[HTML]{333333} $\dot{u} : uv$} & {\color[HTML]{333333} -0.68} & {\color[HTML]{333333} -0.321} & {\color[HTML]{333333} -0.73} & {\color[HTML]{333333} -0.67} \\
\rowcolor[HTML]{FFCE93} 
{\color[HTML]{333333} $\dot{v} : uv$} & {\color[HTML]{333333} 0.82} & {\color[HTML]{333333} 0.71} & {\color[HTML]{333333} 0.78} & {\color[HTML]{333333} 0.73} \\
$\dot{u} : u^2$ & 0 & 0 & 0.00 & 0.00 \\
$\dot{v} : u^2$ & 0 & 0 & 0.00 & 0.00 \\
$\dot{u} : v^2$ & 0 & 0 & 0.00 & 0.00 \\
$\dot{v} : v^2$ & 0 & 0 & 0.00 & 0.00 \\ \hline
\end{tabular}
\quad
\begin{tabular}{|lll|}
\hline
\rowcolor[HTML]{FFFFFF} 
 & ss-SINDy & rh-SINDy \\ \hline
\rowcolor[HTML]{FFFFFF} 
$\dot{u} : 1$ & 0.36 & 0.02 \\
\rowcolor[HTML]{FFFFFF} 
$\dot{v} : 1$ & 0.17 & -0.05 \\
\rowcolor[HTML]{FFCE93} 
{\color[HTML]{333333} $\dot{u} : u$} & {\color[HTML]{333333} 1.00} & {\color[HTML]{333333} 3.38} \\
$\dot{v} : u$ & 0.13 & 0.00 \\
$\dot{u} : v$ & 0.27 & 0.03 \\
\rowcolor[HTML]{FFCE93} 
{\color[HTML]{000000} $\dot{v} : v$} & {\color[HTML]{000000} 1.00} & {\color[HTML]{000000} 1.03} \\
\rowcolor[HTML]{FFCE93} 
{\color[HTML]{000000} $\dot{u} : uv$} & {\color[HTML]{000000} 1.00} & {\color[HTML]{000000} 1.17} \\
\rowcolor[HTML]{FFCE93} 
{\color[HTML]{000000} $\dot{v} : uv$} & {\color[HTML]{000000} 1.00} & {\color[HTML]{000000} 1.17} \\
$\dot{u} : u^2$ & 0.09 & 0.02 \\
$\dot{v} : u^2$ & 0.01 & 0.04 \\
$\dot{u} : v^2$ & 0.19 & 0.04 \\
$\dot{v} : v^2$ & 0.03 & -0.02 \\ \hline
\end{tabular}
\caption{(Left) Posterior modes of SINDy coefficients for the Lotka-Volterra model. (Right) Corresponding inclusion probabilities and pseudo-probabilities.}
\label{table:lotka_volterra}
\end{table}

\begin{figure}[t]
\includegraphics[width=0.95\textwidth]{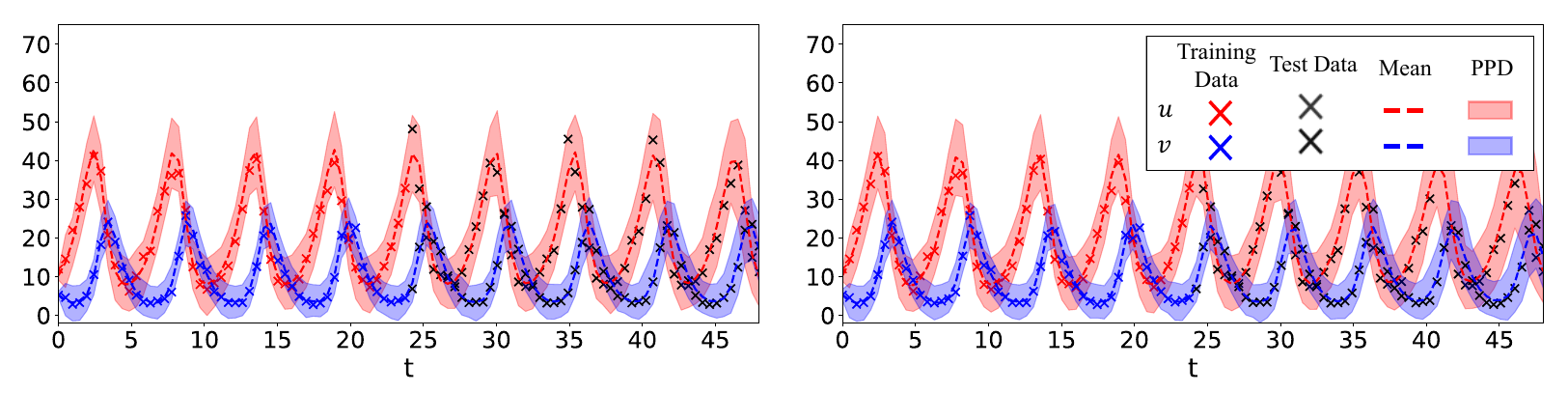}
\caption{Forecasting using ssh-SINDy (left) and rh-SINDy (right) for the Lotka-Volterra model.
  We train using samples the time interval $[0,24]$ (red and blue crosses) and test on samples over the time interval $(24,48]$ (black crosses). The mean (dashed lines) and 90\% credibility intervals (dashed areas) of the PPDs are plotted for the entire time interval.
}
\label{fig:futurestate}
\end{figure}

Finally, we demonstrate how the UQ-SINDy framework can be used for forecasting. To do this, we first simulate noisy data over the time interval $(24,48]$ (black crosses) and use this as our test set (see Figure~\ref{fig:futurestate}). We then compute the PPD over the entire time interval $[0,48]$ by sampling from~\eqref{eq:ppc}, and plot the mean and 90\% credibility interval of this distribution. We find that the mean of the PPD is very close in value to the true values in the test set. Further, we note that some samples in the test set lie near the bounds of the credibility intervals. This shows that our credibility bounds are tight and accurately capture the uncertainty due to measurement noise.

\subsubsection{Nonlinear oscillator and model indeterminacy}\label{subsec:nonlinear_oscillator}

\begin{figure}[t]
\includegraphics[width=0.95\textwidth]{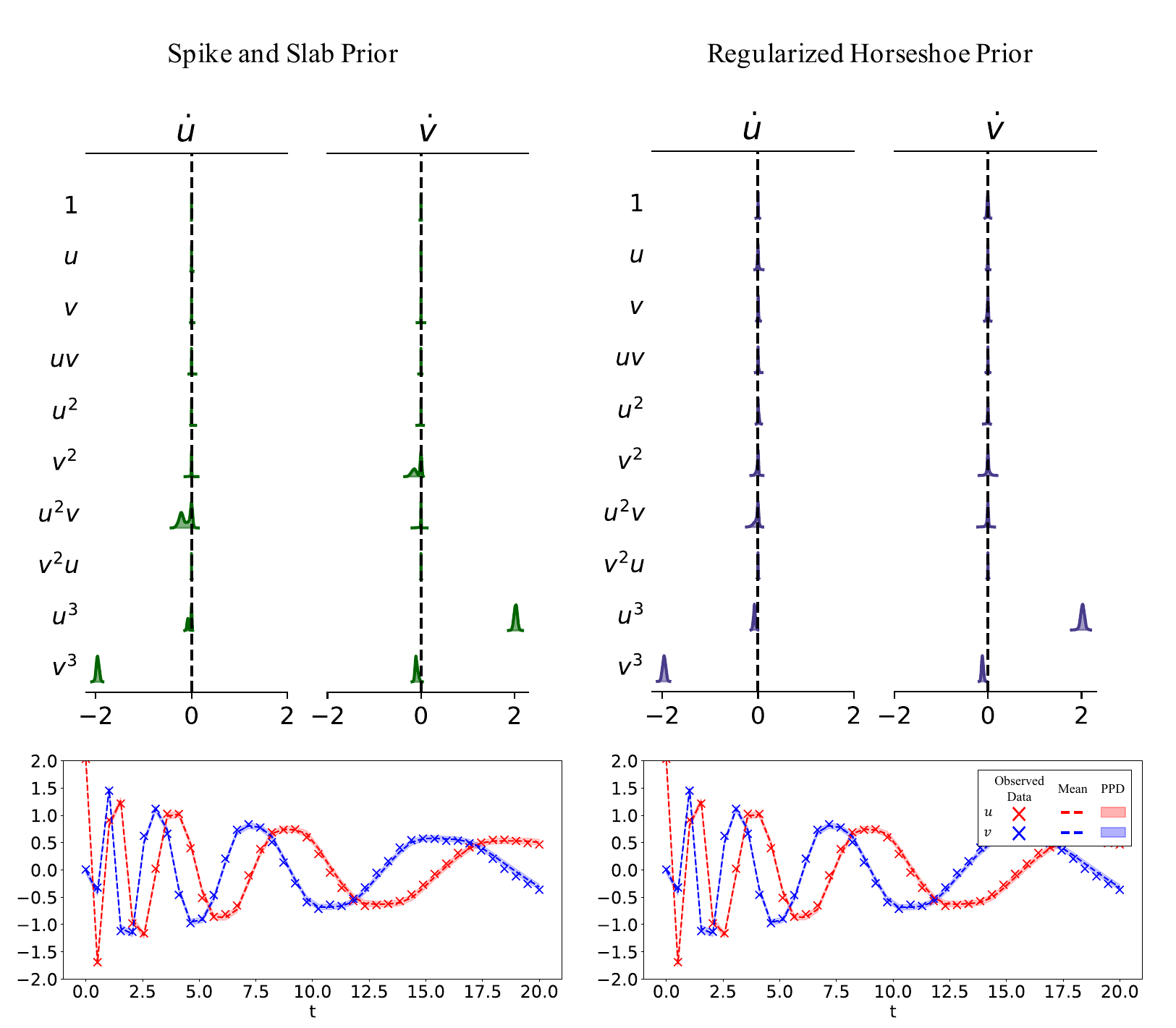}
\vspace*{-.2in}
\caption{UQ-SINDy applied to a synthetic nonlinear oscillator system with normal noise. (Above) Marginal ss-SINDy and rh-SINDy posterior distributions.
(Below) Observed (crosses) and predicted time series together with the corresponding PPD means (dashed lines) and 90\% credibility intervals (shaded areas).
}
\label{fig:nonlinear_oscillator}
\end{figure}

As a second example, we consider the damped nonlinear oscillator model of the form
\begin{align*}
  \dot{u} = \alpha u^3 + \beta v^3 \\
  \dot{v} = \gamma v^3 + \delta u^3.
\end{align*}
Following~\cite{raissi2018hidden}, we use the values $\alpha = -0.1, \beta = -2, \gamma=2, \delta = -0.1$ and the initial conditions $[u_0, v_0] = [2,0]$. Data is generated by sampling this model over the interval $t \in [0, 20]$ with a sampling period of $\Delta t = 0.2$, and adding normally distributed observation noise with distribution $\mathcal{N}(0,0.02^2)$. The observed trajectory is shown in Figure~\ref{fig:nonlinear_oscillator}. We use a library of polynomial terms $\bm{\Theta}(u, v) = [1, u, v, u^2, v^2, u v, u^3, v_3, u^2 v, v^2 u]$, resulting in a $10 \times 2$ matrix of SINDy coefficients $\bm{\Xi}$.
Since the observation noise is normally distributed noise we employ the data likelihood in~\eqref{eq:normal_likelihood}. For the noise level and initial condition we employ the priors $\sigma_u, \sigma_v \sim \operatorname{Gamma}(\alpha = 1, \beta = 0.1)$ and $u_0, v_0 \sim \operatorname{Laplace}(\mu = 0, b = 1)$, respectively.

First, we apply SINDy to the data, resulting in the estimated SINDy coefficients presented in Table~\ref{table:nonlinear_oscillator}.
It can be seen that SINDy does not identify the relevant terms in the model or correctly estimate the values of the model parameters. In fact, none of the terms in the SINDy model are zero. This example is particularly challenging for SINDy because of the sparse data sampling, the size of the library, and the large range of magnitudes of the nonzero coefficients (specifically, note that $|\alpha|$ and $|\delta|$ are much smaller than $|\beta|$ and $|\gamma|$).
Next, we apply ssh-SINDy and rh-SINDy to this data.
The posterior modes are shown in Table~\ref{table:nonlinear_oscillator}.

We present the marginal posterior distributions of the SINDy coefficients in Figure~\ref{fig:nonlinear_oscillator}. It can be seen that rh-SINDy correctly identifies the governing equation; specifically, we see that the marginal posterior distribution of the SINDy coefficients for the terms in the equation are centered away from zero, while the distributions of all other terms are sharply centered at zero.
On the other hand, ss-SINDy identifies the four terms in the governing equation, while also identifying an additional mode corresponding to a model without the $\dot{u} : u^3$ term but with the $\dot{u} : u^2 v$ and $\dot{v} : v^2$ terms.
These results are reflected in Table~\ref{table:nonlinear_oscillator}, for which we show the posterior modes of the SINDy coefficients and the corresponding inclusion probabilities and pseudo-probabilities. For rh-SINDy, the four nonzero terms are clearly identified with modes close to the true values and inclusion pseudo-probabilities for the four terms close to one. For ss-SINDy, three of the terms are clearly identified with an inclusion probability close to one, while the terms $\dot{u} : u^3$, $\dot{u} : u^2 v$ and $\dot{v} : v^2$ have inclusion probabilities of $0.5$, $0.7$, and $0.5$, respectively.

\begin{table}[t]
\centering
\begin{tabular}{|lllll|}
\hline
 & TRUE & SINDy & ss-SINDy & rh-SINDy \\ \hline
$\dot{u} : 1$ & 0 & 0.46 & 0.00 & 0.00 \\
$\dot{v} : 1$ & 0 & 0.06 & 0.00 & 0.00 \\
\rowcolor[HTML]{FFFFFF} 
{\color[HTML]{000000} $\dot{u} : u$} & {\color[HTML]{000000} 0} & {\color[HTML]{000000} 0.54} & {\color[HTML]{000000} 0.00} & {\color[HTML]{000000} 0.00} \\
\rowcolor[HTML]{FFFFFF} 
$\dot{v} : u$ & 0 & -0.62 & 0.00 & 0.00 \\
\rowcolor[HTML]{FFFFFF} 
$\dot{u} : v$ & 0 & 0.81 & 0.00 & 0.00 \\
\rowcolor[HTML]{FFFFFF} 
{\color[HTML]{333333} $\dot{v} : v$} & {\color[HTML]{333333} 0} & {\color[HTML]{333333} -0.09} & {\color[HTML]{333333} 0.00} & {\color[HTML]{333333} 0.00} \\
\rowcolor[HTML]{FFFFFF} 
{\color[HTML]{333333} $\dot{u} : uv$} & {\color[HTML]{333333} 0} & {\color[HTML]{333333} -0.45} & {\color[HTML]{333333} 0.00} & {\color[HTML]{333333} 0.00} \\
\rowcolor[HTML]{FFFFFF} 
{\color[HTML]{333333} $\dot{v} : uv$} & {\color[HTML]{333333} 0} & {\color[HTML]{333333} -0.14} & {\color[HTML]{333333} 0.00} & {\color[HTML]{333333} 0.00} \\
$\dot{u} : u^2$ & 0 & -1.82 & 0.00 & 0.00 \\
$\dot{v} : u^2$ & 0 & -0.38 & 0.00 & 0.00 \\
$\dot{u} : v^2$ & 0 & 0.43 & 0.00 & 0.00 \\
$\dot{v} : v^2$ & 0 & 0.34 & 0.00 & 0.00 \\
{\color[HTML]{333333} $\dot{u} : u^2 v$} & 0 & 0.39 & 0.00 & 0.00 \\
$\dot{v} : u^2 v$ & 0 & 0.15 & 0.00 & 0.00 \\
{\color[HTML]{333333} $\dot{u} : v^2 u$} & 0 & 1.37 & 0.00 & 0.00 \\
{\color[HTML]{333333} $\dot{v} : v^2 u$} & 0 & -0.22 & 0.00 & 0.00 \\
\rowcolor[HTML]{FFCE93} 
{\color[HTML]{333333} $\dot{u} : u^3$} & -0.1 & -1.41 & 0.00 & -0.08 \\
\rowcolor[HTML]{FFCE93} 
$\dot{v} : u^3$ & 2 & -0.53 & 2.04 & 2.02 \\
\rowcolor[HTML]{FFCE93} 
{\color[HTML]{333333} $\dot{u} : v^3$} & -2 & 0.02 & -1.96 & -1.96 \\
\rowcolor[HTML]{FFCE93} 
{\color[HTML]{333333} $\dot{v} : v^3$} & -0.1 & -0.15 & -0.11 & -0.12 \\ \hline
\end{tabular}
\quad
\begin{tabular}{|lll|}
\hline
 & ss-SINDy & rh-SINDy \\ \hline
$\dot{u} : 1$ & 0.00 & -0.04 \\
$\dot{v} : 1$ & 0.01 & -0.19 \\
\rowcolor[HTML]{FFFFFF} 
{\color[HTML]{000000} $\dot{u} : u$} & {\color[HTML]{000000} 0.10} & {\color[HTML]{000000} -0.02} \\
\rowcolor[HTML]{FFFFFF} 
$\dot{v} : u$ & 0.00 & -0.07 \\
\rowcolor[HTML]{FFFFFF} 
$\dot{u} : v$ & 0.05 & 0.03 \\
\rowcolor[HTML]{FFFFFF} 
{\color[HTML]{333333} $\dot{v} : v$} & {\color[HTML]{333333} 0.08} & {\color[HTML]{333333} -0.03} \\
\rowcolor[HTML]{FFFFFF} 
{\color[HTML]{333333} $\dot{u} : uv$} & {\color[HTML]{333333} 0.08} & {\color[HTML]{333333} -0.02} \\
\rowcolor[HTML]{FFFFFF} 
{\color[HTML]{333333} $\dot{v} : uv$} & {\color[HTML]{333333} 0.06} & {\color[HTML]{333333} 0.00} \\
$\dot{u} : u^2$ & 0.07 & 0.01 \\
$\dot{v} : u^2$ & 0.07 & 0.00 \\
$\dot{u} : v^2$ & 0.14 & 0.25 \\
$\dot{v} : v^2$ & 0.5 & -0.01 \\
{\color[HTML]{333333} $\dot{u} : u^2 v$} & 0.70 & 0.47 \\
$\dot{v} : u^2 v$ & 0.17 & 0.01 \\
{\color[HTML]{333333} $\dot{u} : v^2 u$} & 0.01 & 0.38 \\
{\color[HTML]{333333} $\dot{v} : v^2 u$} & 0.02 & 0.14 \\
\rowcolor[HTML]{FFCE93} 
{\color[HTML]{333333} $\dot{u} : u^3$} & 0.50 & 1.24 \\
\rowcolor[HTML]{FFCE93} 
$\dot{v} : u^3$ & 1.00 & 1.05 \\
\rowcolor[HTML]{FFCE93} 
{\color[HTML]{333333} $\dot{u} : v^3$} & 1.00 & 0.99 \\
\rowcolor[HTML]{FFCE93} 
{\color[HTML]{333333} $\dot{v} : v^3$} & 1.00 & 0.82 \\ \hline
\end{tabular}
\caption{(Left) Posterior modes of SINDy coefficients for the nonlinear oscillator model. (Right) Corresponding inclusion probabilities and pseudo-probabilities.}
\label{table:nonlinear_oscillator}
\end{table}

\begin{figure}[b]
\includegraphics[width=\textwidth]{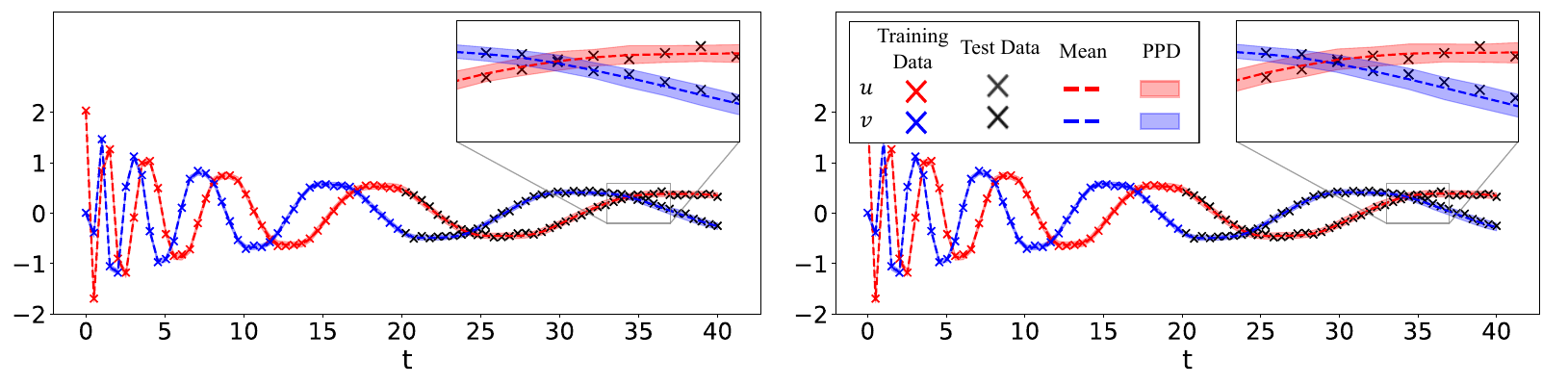}
\vspace*{-.3in}
\caption{Forecasting using ssh-SINDy (left) and rh-SINDy (right) for the nonlinear oscillator model.
  We train using samples from the Lotka-Volterra model over the time interval $[0,20]$ (red and blue crosses) and test on samples over the time interval $(20,40]$ (black crosses). The mean (dashed lines) and 90\% credibility intervals (dashed areas) of the PPDs are plotted for the entire time interval.
}
\label{fig:nonlinear_oscillator_future_state}
\end{figure}

In Figure~\ref{fig:nonlinear_oscillator}, we present the mean and 90\% credibility intervals of the PPDs of the reconstruction of the system's states, together with the training data.
Similarly, in Figure~\ref{fig:nonlinear_oscillator_future_state} we present the 90\% credibility intervals of the PPDs of future state forecasting for testing data over the time interval $(20,40]$.
Both rh-SINDy and ss-SINDy lead to similar credibility intervals for both reconstruction and forecasting.
We note that the range of predicted model states is much narrower than for the Lotka-Volterra model, which is expected due to the lower noise level present in these measurements. We also emphasize that these PPDs are much tighter than those presented in~\cite{yang2020bayesian} for this test case, even though we train rh-SINDy and ss-SINDy with substantially less data than in that work.
Furthermore, it can be seen that the test data lies within the 90\% credibility intervals of the PPDs of each state. Although some of the draws from the ss-SINDy PPD contain terms not in the model, the credibility intervals for both ss-SINDy and rh-SINDy are similar. This suggests that the ambiguity identified by the spike and slab prior is due to model indeterminacy inherent in this data set.

\begin{figure}[t]
\centering
\includegraphics[width=0.8\textwidth]{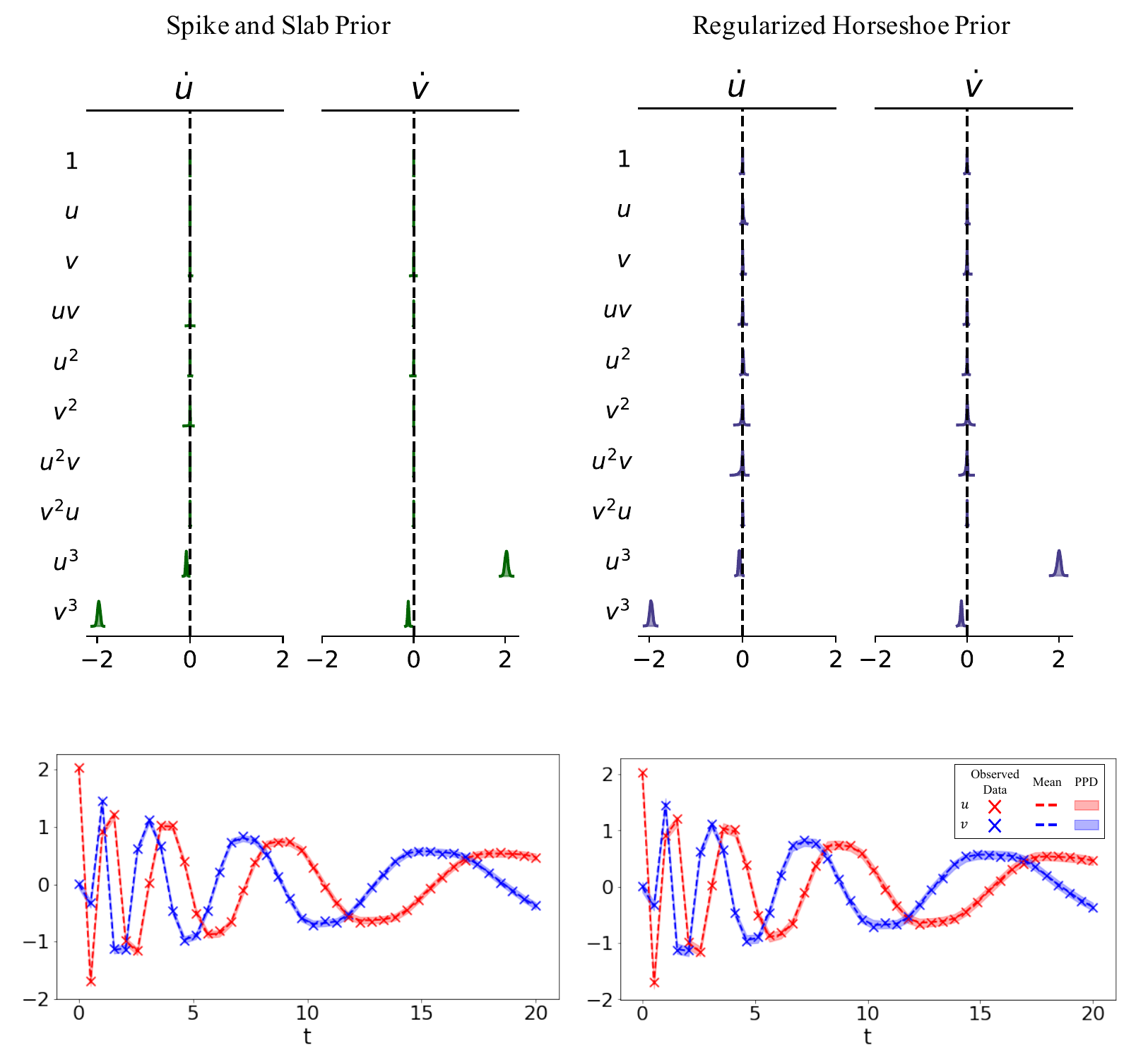}
\caption{UQ-SINDy, with scaled priors for the terms $\dot{u} : u^3$ and $\dot{v} : v^3$, applied to a synthetic nonlinear oscillator system with normal noise. (Above) Marginal ss-SINDy and rh-SINDy posterior distributions.
  (Below) Observed (crosses) and predicted time series together with the corresponding PPD means (dashed lines) and 90\% credibility intervals (shaded areas).
}
\label{fig:nonlinear_oscillator_scaled}
\end{figure}

This indeterminacy can be attributed to the range of values spanned by the coefficients in the governing equation. In particular, the coefficients of $\dot{u} : u^3$ and $\dot{v} : v^3$ are an order of magnitude smaller than the coefficients of $\dot{u} : v^3$ and $\dot{u} : u^3$.
To further investigate this indeterminacy, we re-applied ss-SINDy and rh-SINDy with $\alpha_{i, j} = 0.1$ for the terms $\dot{u} : u^3$ and $\dot{v} : v^3$.
This scaling of the prior incorporates the knowledge that these two terms have coefficients of magnitude $O(0.1)$.
The resulting marginal posterior distributions, presented in Figure~\ref{fig:nonlinear_oscillator_scaled}, show that this scaling of the prior removes this ambiguity.
The corresponding posterior modes and inclusion probabilities and pseudo-probabilities are presented in Table~\ref{table:nonlinear_oscillator_scaled}.

\begin{table}[t]
\centering
\begin{tabular}{|llll|}
\hline
 & TRUE & ss-SINDy & rh-SINDy \\ \hline
$\dot{u} : 1$ & 0 & 0.00 & 0.00 \\
$\dot{v} : 1$ & 0 & 0.00 & 0.00 \\
$\dot{u} : u$\} & 0 & 0.00 & 0.00 \\
$\dot{v} : u$ & 0 & 0.00 & 0.00 \\
$\dot{u} : v$ & 0 & 0.00 & 0.00 \\
$\dot{v} : v$ & 0 & 0.00 & 0.00 \\
$\dot{u} : uv$ & 0 & 0.00 & 0.00 \\
$\dot{v} : uv$ & 0 & 0.00 & 0.00 \\
$\dot{u} : u^2$ & 0 & 0.00 & 0.00 \\
$\dot{v} : u^2$ & 0 & 0.00 & 0.00 \\
$\dot{u} : v^2$ & 0 & 0.00 & 0.00 \\
$\dot{v} : v^2$ & 0 & 0.00 & 0.00 \\
$\dot{u} : u^2 v$ & 0 & 0.00 & 0.00 \\
$\dot{v} : u^2 v$ & 0 & 0.00 & 0.00 \\
$\dot{u} : v^2 u$ & 0 & 0.00 & 0.00 \\
$\dot{v} : v^2 u$ & 0 & 0.00 & 0.00 \\
\rowcolor[HTML]{FFCE93} 
$\dot{u} : u^3$ & -0.1 & -0.08 & -0.07 \\
\rowcolor[HTML]{FFCE93} 
$\dot{v} : u^3$ & 2 & 2.03 & 2.01 \\
\rowcolor[HTML]{FFCE93} 
$\dot{u} : v^3$ & -2 & -1.97 & -1.97 \\
\rowcolor[HTML]{FFCE93} 
$\dot{v} : v^3$ & -0.1 & -0.12 & -0.12 \\
\hline
\end{tabular}
\quad
\begin{tabular}{|lll|}
\hline
 & ss-SINDy & rh-SINDy \\ \hline
$\dot{u} : 1$ & 0.00 & 0.01 \\
$\dot{v} : 1$ & 0.00 & -0.01 \\
$\dot{u} : u$ & 0.00 & 0.02 \\
$\dot{v} : u$ & 0.00 & 0.03 \\
$\dot{u} : v$ & 0.03 & 0.01 \\
$\dot{v} : v$ & 0.11 & -0.01 \\
$\dot{u} : uv$ & 0.10 & -0.03 \\
$\dot{v} : uv$ & 0.10 & 0.01 \\
$\dot{u} : u^2$ & 0.08 & 0.08 \\
$\dot{v} : u^2$ & 0.05 & 0.00 \\
$\dot{u} : v^2$ & 0.17 & 0.11 \\
$\dot{v} : v^2$ & 0.00 & 0.014 \\
$\dot{u} : u^2 v$ & 0.00 & 0.01 \\
$\dot{v} : u^2 v$ & 0.00 & 0.03 \\
$\dot{u} : v^2 u$ & 0.03 & -0.56 \\
$\dot{v} : v^2 u$ & 0.01 & -2.53 \\
\rowcolor[HTML]{FFCE93} 
$\dot{u} : u^3$ & 1.00 & 1.23 \\
\rowcolor[HTML]{FFCE93} 
$\dot{v} : u^3$ & 1.00 & 1.05 \\
\rowcolor[HTML]{FFCE93} 
$\dot{u} : v^3$ & 1.00 & 1.00 \\
\rowcolor[HTML]{FFCE93} 
$\dot{v} : v^3$ & 1.00 & 0.93 \\
\hline
\end{tabular}
\caption{(Left) Posterior modes of SINDy coefficients, with scaled priors for the terms $\dot{u} : u^3$ and $\dot{v} : v^3$, for the nonlinear oscillator model.
  (Right) Corresponding inclusion probabilities and pseudo-probabilities.}
\label{table:nonlinear_oscillator_scaled}
\end{table}

\subsubsection{Lynx-hare population model}\label{subsec:lynx_hare}

\begin{figure}[t]
\includegraphics[width=\textwidth]{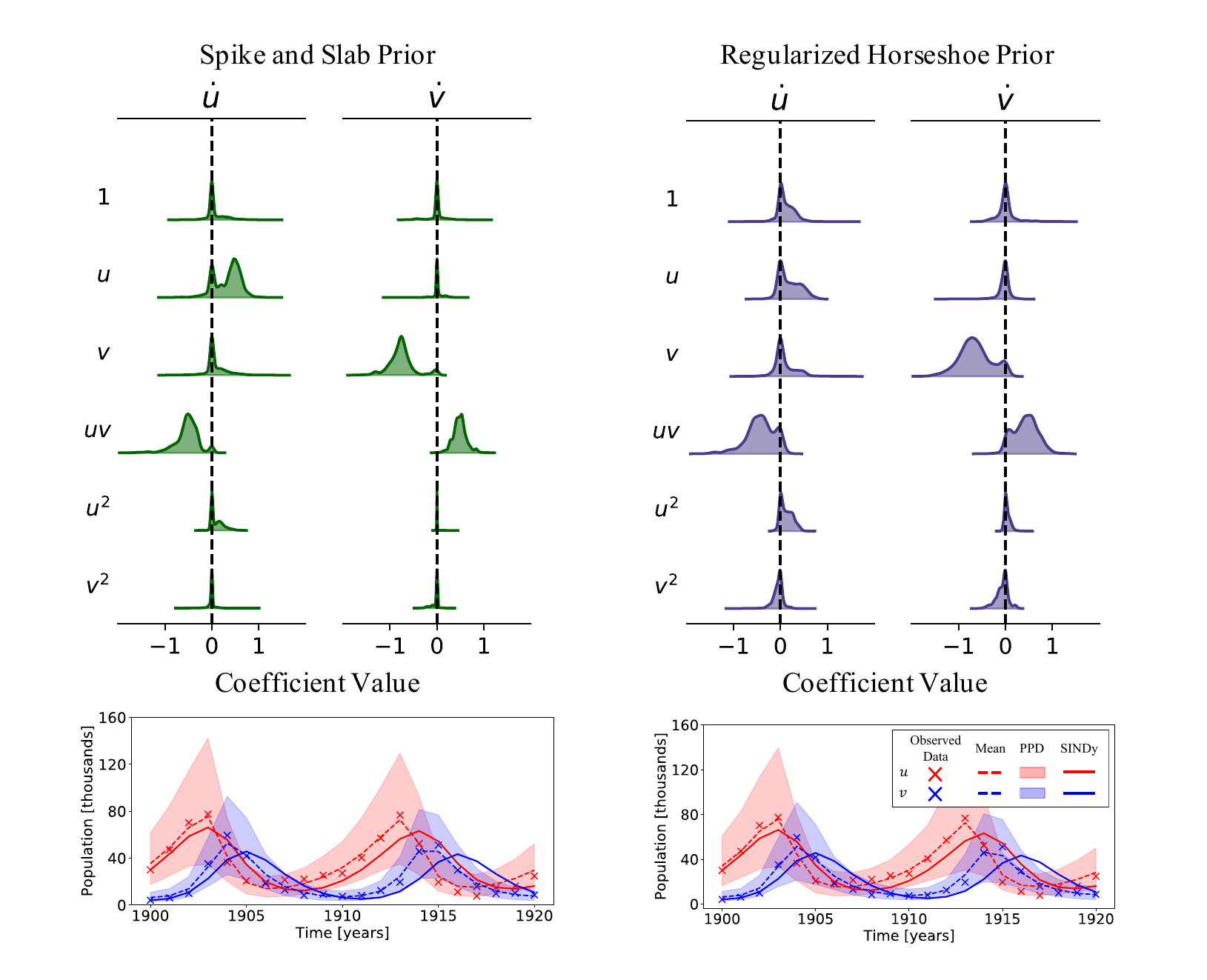}
\vspace*{-.5in}
\caption{UQ-SINDy applied to the lynx-hare population data.
  (Above) Marginal ss-SINDy and rh-SINDy posterior distributions.
  (Below) Observed (crosses) and predicted time series together with the corresponding PPD means (dashed lines) and 90\% credibility intervals (shaded areas).
  SINDy predictions presented as continuous lines.
}
\label{fig:predator_prey}
\end{figure}

\begin{table}[t]
\centering
\begin{tabular}{|lllll|}
\hline
 & Param. est. & SINDy & ss-SINDy & rh-SINDy \\ \hline
$\dot{u} : 1$ & 0 & 0 & 0.00 & 0.01 \\
$\dot{v} : 1$ & 0 & 0 & 0.00 & 0.00 \\
\rowcolor[HTML]{FFCE93} 
{\color[HTML]{333333} $\dot{u} : u$} & {\color[HTML]{333333} 0.55} & {\color[HTML]{333333} 0.48} & {\color[HTML]{333333} 0.47} & {\color[HTML]{333333} 0.00} \\
$\dot{v} : u$ & 0 & 0 & 0.00 & -0.01 \\
$\dot{u} : v$ & 0 & -0.143 & 0.00 & 0.00 \\
\rowcolor[HTML]{FFCE93} 
{\color[HTML]{333333} $\dot{v} : v$} & {\color[HTML]{333333} -0.84} & {\color[HTML]{333333} -0.71} & {\color[HTML]{333333} -0.76} & {\color[HTML]{333333} -0.7} \\
\rowcolor[HTML]{FFCE93} 
{\color[HTML]{333333} $\dot{u} : uv$} & {\color[HTML]{333333} -0.455} & {\color[HTML]{333333} -0.36} & {\color[HTML]{333333} -0.51} & {\color[HTML]{333333} -0.42} \\
\rowcolor[HTML]{FFCE93} 
{\color[HTML]{333333} $\dot{v} : uv$} & {\color[HTML]{333333} 0.5433} & {\color[HTML]{333333} 0.42} & {\color[HTML]{333333} 0.52} & {\color[HTML]{333333} 0.52} \\
$\dot{u} : u^2$ & 0 & 0 & 0.00 & 0.01 \\
$\dot{v} : u^2$ & 0 & 0 & 0.00 & 0.00 \\
$\dot{u} : v^2$ & 0 & 0 & 0.00 & 0.00 \\
$\dot{v} : v^2$ & 0 & 0 & 0.00 & -0.01 \\ \hline
\end{tabular}
\quad
\begin{tabular}{|lll|}
\hline
 & ss-SINDy & rh-SINDy \\ \hline
$\dot{u} : 1$ & 0.47 & 0.04 \\
$\dot{v} : 1$ & 0.35 & 0.00 \\
\rowcolor[HTML]{FFCE93} 
{\color[HTML]{000000} $\dot{u} : u$} & {\color[HTML]{000000} 0.85} & {\color[HTML]{000000} 0.01} \\
$\dot{v} : u$ & 0.34 & 0.02 \\
$\dot{u} : v$ & 0.54 & 0.00 \\
\rowcolor[HTML]{FFCE93} 
{\color[HTML]{333333} $\dot{v} : v$} & {\color[HTML]{333333} 0.99} & {\color[HTML]{333333} 0.73} \\
\rowcolor[HTML]{FFCE93} 
{\color[HTML]{333333} $\dot{u} : uv$} & {\color[HTML]{333333} 0.96} & {\color[HTML]{333333} 0.78} \\
\rowcolor[HTML]{FFCE93} 
{\color[HTML]{333333} $\dot{v} : uv$} & {\color[HTML]{333333} 1} & {\color[HTML]{333333} 2.01} \\
$\dot{u} : u^2$ & 0.581 & 0.03 \\
$\dot{v} : u^2$ & 0.08 & 0.02 \\
$\dot{u} : v^2$ & 0.31 & -0.50 \\
$\dot{v} : v^2$ & 0.35 & -0.06 \\ \hline
\end{tabular}
\caption{(Left) Posterior modes of SINDy coefficients for the lynx-hare data. (Right) Corresponding inclusion probabilities and pseudo-probabilities.}
\label{table:predator_prey}
\end{table}

As a final example, we apply ss-SINDy and rh-SINDy as described in Section~\ref{subsec:lotka_volterra} to model the population dynamics of two species in Canada. In particular, we consider data consisting of measurements by the Hudson Bay Company of lynx and hare pelts between 1900 and 1920~\cite{carpenter2018predator,hewitt1921conservation} (see Figure~\ref{fig:predator_prey}). The number of pelts for these two species is thought to be proportional to the true populations. Hares are a herbivorous relative of the rabbit, while the lynx is a type of wildcat whose diet depends heavily on hares. This predator-prey interdependence between the two species has been shown to be well characterized to first-order by the Lotka-Volterra model~\eqref{eq:lotka_volterra}, where $u$ and $v$ correspond to the populations of hares and lynx, respectively.

Figure~\ref{fig:predator_prey} presents the number of pelts recorded yearly for these two species over 21 years. Modeling this data with SINDy is particularly challenging because we have relatively few samples that cover only two cycles. In addition, factors such as the weather and the consistency of trapping between years adds uncertainty to the measurements.
Here we compare the performance of ss-SINDY, and rh-SINDY for model discovery under uncertainty. The SINDy library, as in the Lotka-Volterra example, contains all constant, linear and quadratic terms. In addition, as a preprocessing step we normalize the data as described in Section~\ref{subsec:lotka_volterra}

The marginal posterior distributions computed using ss-SINDy and rh-SINDy are presented in Figure~\ref{fig:predator_prey}.
The posterior modes and inclusion probabilities and pseudo-probabilities are presented in Table~\ref{table:predator_prey}, together with maximum likelihood estimates of the coefficients of the Lotka-Volterra model for the lynx-hare data~\cite{carpenter2018predator}, and estimates computed using the original SINDy algorithm.
It can be seen that for ss-SINDy the distinct nonzero peaks corresponding to the terms in~\eqref{eq:lotka_volterra}. The likelihood of these four terms belonging to the model are very high. We additionally see a small peak near zero for $\dot{u} : u$. This term is highly correlated with a nonzero constant term. We see a similar but more pronounced peak for rh-SINDy.
Table~\ref{table:predator_prey} shows that ss-SINDy correctly identifies the Lotka-Volterra model and assigns high inclusion probabilities to the four terms in such a model.
On the other hand, rh-SINDy identifies three of the four terms correctly. 
Furthermore, it can be seen that SINDy fails to identify the Lotka-Volterra model, and that the posterior modes for ss-SINDy and rh-SINDy are closer to the maximum likelihood estimates than the SINDy estimates.

Last, in Figure~\ref{fig:predator_prey}, we present the mean and 90\% credibility intervals of the PPDs of the time series reconstruction. We note that all data lie within these credibility bounds.
The SINDy reconstructions, on the other hand, appear to deviate from the time series for later times.

\section{Conclusions and future work}\label{sec:conclusions}

In this work we proposed UQ-SINDy, a new uncertainty quantification framework for identifying governing ODEs directly from noisy and sparse time series data.
We leverage advances in model discovery for dynamical systems and sparsity promoting Bayesian inference to identify a sparse set of SINDy library functions that best explain the observed data, and to quantifying the uncertainty in the SINDy coefficients due to measurement noise and the probability of inclusion of each term in the SINDy library into the final model. We have applied UQ-SINDy to two synthetic examples and one real-world example of lynx-hare population data. By utilizing the spike-and-slab and regularized horseshoe priors, UQ-SINDy yields posterior distributions of SINDy coefficients with truly sparse draws, and thus results in truly sparse probabilistic model discovery; in contrast, the use of the Laplace prior does not lead to sparse model discovery. We observe that the proposed approach is robust against observation noise and can accommodate sparse samples and small data sets.

Going forward, one of the primary limitations of this method is its scalability to very large SINDy libraries. This is primarily due to the computational cost of sampling high-dimensional posterior distributions using MCMC. One remedy for this is to use variational inference, which matches classes of distributions to the posterior distribution by maximizing a lower bound on the marginal likelihood of the data. This method has been particuarly effective for high dimensional models, most notably neural networks, with comparable accuracy to sampling-based methods. 
Furthermore, in this work we are primarily focused on situations in which the coordinates that induce a sparse representation are known. However, in general this ``effective'' set of coordinates may be unknown. Recent work merges SINDy together with neural network architectures in order to simultaneously learn parsimonous governing equations and the associated sparsity-inducing coordinate transformation~\cite{champion2019data}. Extending UQ-SINDy to this coordinate discovery framework could greatly improve the robustness of the learning process under uncertainty and the quality of the resulting forecasts.

\section*{Acknowledgements}

This research was supported by Laboratory Directed Research and Development Program and Mathematics for Artificial Reasoning for Scientific Discovery investment at the Pacific Northwest National Laboratory, a multiprogram national laboratory operated by Battelle for the U.S. Department of Energy under Contract DE-AC05- 76RLO1830.

\bibliographystyle{plain}
\bibliography{bayesian_sindy_citations}

\appendix

\end{document}